\documentclass{article}
\usepackage{amsmath,amssymb,amsthm}
\usepackage{graphics}

\newtheorem{theorem}{Theorem}
\newtheorem{lemma}{Lemma}

\newcommand{\conv}{\mathop{\rm conv}\nolimits}
\newcommand{\dist}{\mathop{\rm dist}\nolimits}

\newtheorem{cor}{Corollary}

\title {The one-sided kissing number in four dimensions}

\author {Oleg R. Musin \thanks{Institute for Math. Study of Complex Systems, Moscow State University, Moscow, Russia omusin@mail.ru}}

\begin{document}
\date{}
\maketitle

\begin{abstract}
Let $H$ be a closed half-space of $n$-dimensional Euclidean space. Suppose $S$ is a unit sphere in $H$ that touches the supporting hyperplane of $H$. The one-sided kissing number $B(n)$ is the maximal number of unit nonoverlapping spheres in $H$ that can touch $S$. Clearly, $B(2)=4$. It was proved that
$B(3)=9.$ Recently, K. Bezdek proved that $B(4)=18$ or 19, and conjectured that $B(4)=18.$ We present a proof of this conjecture.
\end {abstract}

\section {Introduction}

%Recently, K\'aroly Bezdek \cite {B1,KB} introduced the one-sided kissing problem in $n$-dimensions.

Let $H$ be a closed half-space of ${\bf R}^{n}$. Suppose $S$ is a unit sphere in $H$ that touches the supporting hyperplane of $H$. The one-sided kissing number $B(n)$ is the maximal number of unit nonoverlapping spheres in $H$ that can touch $S$. 

The problem of finding $B(3)$ was raised by L. Fejes T\'oth and H. Sachs in 1976 \cite{FS}. 
%Research problem 17, Periodica Math. Hungar. 7 (1976), 125-127. L. Fejes T\'oth
%came to the problem through the conjecture that any 10-neighbour packing
%of equal balls in R^3 has positive density, while Horst Sachs arrived to
%the problem studying the chromatic number of finite ballpackings. They did
%not introduced a term for the number, and by the nature of the Research
%Problem Secion of the journal they emphasised the simplest unsolved
%problem and did not mention the analogous question in higher dimensions.
%In their joint paper On k^+-neighbour packings and one-sided Hadwiger
%configurations, Contributions to Algebra and Geometry, 4 (2003), 493-498,
K. Bezdek and P. Brass \cite{BB} studied the problem in more
general setting and they introduced the term ``one-sided Hadwiger
number", which in the case of a ball is the same as the one-sided kissing
number. The term ``one-sided kissing number" have been introduced by K. Bezdek \cite{B1, KB}.

\begin{center}
\begin{picture}(320,100)(-160,-100)
%Fig.1 

\put(-150,-80){\line(1,0){300}}

\put(0,-60){\circle{40}}
\put(-7,-65){$S$}
\put(87,-40){$H$}

\put(-40,-60){\circle{40}}
\put(40,-60){\circle{40}}

\put(-20,-25){\circle{41}}
\put(20,-25){\circle{41}}

\put(20,-60){\circle*{4}}
\put(-20,-60){\circle*{4}}
\put(-10,-43){\circle*{4}}
\put(10,-43){\circle*{4}}

\put(-13,-98) {Fig. 1}

\end{picture}
\end{center}
%--------- end of Fig. 1------------------------------------------

The one-sided kissing problem looks more complicated than the classical kissing number problem. Actually, that is a constrained kissing problem. Denote by $k(n)$ the kissing number in $n-$dimensions, i.e. $k(n)$ is the highest number of equal nonoverlapping spheres in ${\bf R}^n$ that can touch another sphere of the same size. Currently, 
$k(n)$ known for $n=2,3,4,8,24$: $k(2)=6$, $k(3)=12$ \cite{SvdW2}, $k(4)=24$ \cite{Mus,Mus2}, $k(8)=240, \; k(24)=196,560$ \cite{Lev2,OdS}.  
There are strong relations between $k(n)$ and $B(n)$ (see the Conclusions and Theorem 1 in this paper).

If $M$ unit spheres kiss (touch) the unit sphere in $H\subset{\bf R}^n$, then the set of kissing points 
is an arrangement on the hemisphere $S_+$ of $S$ such that the (Euclidean) distance between any two points is at least 1. So the one-sided kissing number problem can be stated in another way: How many points can be placed on the surface of $S_+$  so that the angular separation between any two points is at least $\pi/3$?

%In other words,$B(n)=A(n,\pi/3,\pi/2.$ 
Clearly, $B(2)=4$ (see Fig. 1).  It have been proved by G. Fejes T\'oth \cite{GFT} in 1981 that $B(3)=9$  (see also  H. Sachs \cite{Sachs} and   A. Bezdek \& K. Bezdek \cite{AKB} for other proofs). 
G. Kert\'esz (1994) proved that the maximal one-sided kissing arrangement is unique (up to isometry) on the 2-hemisphere $S_+$.

Note that $B(4)\ge 18.$ Indeed, let $A=1/\sqrt{2}$,
$$ p_1=(A,0,0,A), \quad p_2=(-A,0,0,A),$$
$$p_{\{3,4\}}=(0,\pm A,0,A),\quad 
p_{\{5,6\}}=(0,0,\pm A,A), $$
$$ p_{\{7,...,10\}}=(\pm A,\pm A,0,0),\quad  p_{\{11,...,14\}}=(\pm A,0,\pm A,0),$$
$$p_{\{15,...,18\}}=(0,\pm A,\pm A,0).$$
Then these 18 points $\{p_i\}$ lie on the closed upper hemisphere of the unit sphere with minimal angular separation 
$\pi/3$.\footnote {Actually, the set $\{p_i\}$ is the part of the vertices of the ``24-cell" that belong to $S_+$.}

L. Szab\'o (1991) using the Odlyzko - Sloane bound $k(4)\le 25$ \cite{OdS}  has proved that $B(4) \le 20$. 
Recently, K. Bezdek \cite{B1,KB} based on $k(4)=24$ showed that $B(4)\le 19$  and conjectured that $B(4)=18.$  
% and methods similar to \cite{AKB}

We present in this paper a proof that $B(4)=18$ which is based on our extension of the Delsarte method.
% \cite{Mus,Mus13,Mus2}.

\section {On one-sided kissing arrangement}

%Let us consider an one-sided kissing arrangement $P=\{p_i\}\subset S_+\subset {\bf S}^{n-1}$. 
We denote 
by ${\bf S}^{n-1}$ the unit sphere in ${\bf R}^n$, 
$${\bf S}^{n-1}:=\{(x_1,\ldots,x_{n-1},x_n)\in {\bf R}^n: x_1^2+...+x_n^2=1\}$$
Let $S_+$ be the closed hemisphere of ${\bf S}^{n-1}$ with $x_n\ge 0,$  i.e $S_+$ is a spherical cap of center $N$  and angular radius $90^\circ$, where the ``North pole" $N=(0,\ldots,0,1)$.
The boundary of $S_+$ is called the equator and is denoted  by $E$. ($E\cong{\bf S}^{n-2}$.) 

Let us for $p=(x_1,\ldots,x_{n-1},x_n)$ denote by $p'=(x_1,\ldots,x_{n-1},-x_n)$, and let for $P\subset S_+$ denote by
$P'=\{p': p\in P\}$.

Let $\Pi_n:\dot S_+\to E, \; \dot S_+:= S_+\setminus \{N\}$, be the projection  from the pole $N$ which sends the point $p\in \dot S_+$ along its meridian to the equator of the sphere. 
%In other words, $\Pi_n$ is the projection 
% that takes each point $p=(x_1,\ldots,x_{n-1},x_n)$ to $q=\Pi_n(p)=(y_1,\ldots,y_{n-1},0)$. 
$Q=\{q_1,\ldots,q_m\}\subset E$ is said to be the projection of 
$P=\{p_1,\ldots,p_m\}\subset \dot S_+\; $ if $\; q_i=\Pi_n(p_i).$ Then $q_i=q_j$ if and only if $p_i$ and $p_j$ lie in the same meridian.

Throughout this paper we use the function $\omega(\varphi,\alpha,\beta)$ defined by
%Let $ABC$ be a spherical triangle with sides of angular lengths $\phi_C, \phi_A, \phi_B$. From the spherical law of cosines we have $\angle BAC=\omega(\phi_A,\phi_B,\phi_C)$, where
$$\cos{\omega(\varphi,\alpha,\beta)}:=\frac{\cos{\varphi}-\cos{\alpha}\cos{\beta}}{\sin{\alpha}\sin{\beta}}.$$
This function describes the change of distance $\varphi$ between two points on $\dot S_+$ which are $\alpha$ and $\beta$ away from $N$ under the action $\Pi_n$.

For fixed $\varphi\le 90^\circ$ and $0<u\le\alpha\le\beta\le v\le 90^\circ$, the function $\omega(\varphi,\alpha,\beta)$ attains its minimum at $\alpha=\beta=v\; $ or $\; \alpha=u, \; \beta=v.$ If $v<\varphi$, then the minimum of $\omega$ achieves at $\alpha=\beta=v$ (see Theorem 3 in \cite {Mus2}).

A set $P$ in ${\bf S}^{n-1}$ is called $\varphi$-code if the angular separation between any two points from $P$ is at least $\varphi.$ Denote by $A(n,\varphi)$ the maximal size of  $\varphi$-codes in  ${\bf S}^{n-1}$. 
By $A(n,\varphi,\psi)$ denote the maximal size of $\varphi$-codes  in a spherical cap with angular radius $\psi$. So we have $A(n,\varphi)=A(n,\varphi,180^\circ), \; k(n)=A(n,60^\circ),\; B(n)=A(n,60^\circ,90^\circ).$

We say that a set $P\subset {\bf S}^{n-1}$ is a {\em kissing arrangement} in $\Delta\subset {\bf S}^{n-1}$ if $P\subset\Delta$ and $P$ is a $60^\circ$-code.
$P$ is called a {\em one-sided kissing arrangement} if $\Delta=S_+$. 

\medskip

\noindent{\bf Definition.} Let $P\subset S_+$. Then
$$ P_a:= \{p\in P: \;  90^\circ\ge\dist(p,N)\ge 60^\circ\},\quad a(P):=|P_a|,$$
$$ P_b= \{p\in P: \;  \dist(p,N)< 60^\circ\},\quad b(P):=|P_b|.$$
So we have
$$P=P_a\bigcup P_b, \quad |P|=a(P)+b(P).$$

\medskip

\begin{theorem}
Let $P\subset {\bf S}^{n-1}$ be a one-sided kissing arrangement. Then
 $$a(P)+2b(P)\le k(n).$$
\end{theorem}

\begin{proof}
For any distinct points $p,q \in P$ we have  $\dist(p,q)\ge 60^\circ$.
 It is easy to see (Fig. 2) that $\dist(p',q)\ge\dist(p,q)\ge 60^\circ$. Note that $\dist(P_b,P'_b)\ge 60^\circ.$ Therefore,
$\tilde P=P\bigcup P'_b$ is a kissing arrangement in  ${\bf S}^{n-1}$. Thus, $|\tilde P|\le k(n).$ 
\end{proof}

%\normalsize
\begin{center}
\begin{picture}(320,120)(-240,-70)
%Fig.1 

\put(-140,20){\circle*{3}}
\put(-140,-20){\circle*{3}}
\put(-40,20){\circle*{3}}

\thicklines
\put(-200,0){\line(1,0){170}}

\thinlines
\put(-40,20){\line(-1,0){100}}
\put(-40,20){\line(-5,-2){100}}

\put(-138,-36){$p'$}
\put(-138,28){$p$}
\put(-38,28){$q$}
\put(-185,-30){$S_-$}
\put(-185,30){$S_+$}

\put(-118,-55){Fig. 2}

\end{picture}
\end{center}

\begin{theorem} Let $P\subset {\bf S}^{n-1}$ be a one-sided kissing arrangement. Then
$$a(P)\le A(n-1,\arccos(1/\sqrt{3})).$$
\end{theorem}
\begin{proof} 
Let $Q_a=\{q_i\}$ be the projection of $P_a=\{p_i\},\; q_i=\Pi_n(p_i)$ (see Fig. 3), and let
$$\theta_i=\dist(N,p_i),\; \phi_{i,j}=\dist(p_i,p_j), \; \gamma_{i,j}=\dist(q_i,q_j).$$ Then $\gamma_{i,j}=\omega(\phi_{i,j},\theta_i,\theta_j)$.
We have $\phi_{i,j}\ge 60^\circ$, then  $\gamma_{i,j}\ge\omega(60^\circ,\theta_i,\theta_j)$.
%$$\cos{\gamma_{i,j}}\le \frac{1/2-\cos{\theta_i}\cos{\theta_j}}{\sin{\theta_i}\sin{\theta_j}}.$$
Since $60^\circ\le\theta_i\le 90^\circ$ for all $i$, we see that  $$\gamma_{i,j}\ge \lambda=
\omega(60^\circ,60^\circ,90^\circ)=\arccos(1/\sqrt{3}).$$ Then $Q_a$ is a $\lambda$-code in  ${\bf S}^{n-2}$. Thus, 
$a(P)=|P_a|=|Q_a|\le A(n-1,\lambda).$
\end{proof}

%\medskip

%Fig. 2
\begin{center}
\begin{picture}(320,140)(-60,-70)
\put(80,-63){Fig. 3}
%\qbezier(30,-20)(90,140)(150,-20) 
%\qbezier(30,-20)(29,0)(40,20)
%\qbezier(150,-20)(151,0)(140,20)

%\qbezier(78,58)(90,62)(102,58)
\qbezier(67,54)(90,65)(113,54)
\qbezier (67,54) (55,47) (48,37)
\qbezier (113,54) (125,47) (132,37)

\qbezier (48,37) (38,24) (35,11)
\qbezier (35,11) (31,-4) (30,-20)

\qbezier (132,37) (142,24) (145,11)
\qbezier (145,11) (149,-4) (150,-20)

\put(90,60){\circle*{4}}
\put(83,66){$N$}

\put(62,14){\circle*{4}}
\put(121,10){\circle*{4}}

\put(45,-24){\circle*{4}}
\put(135,-24){\circle*{4}}
\put(40,-34){$q_i$}
\put(134,-34){$q_j$}

%\put(170,120){\circle*{4}}

\thicklines
\qbezier(30,-20)(90,-40)(150,-20)
\thinlines

\qbezier(30,-20)(90,0)(150,-20)

\qbezier(90,60)(67,35)(45,-25)
\qbezier(90,60)(113,35)(135,-25)
\qbezier(62,14)(90,15)(121,10)

\put(49,16){$p_i$}
\put(125,10){$p_j$}
\put(85,5){$\phi_{i,j}$}
\put(85,-39){$\gamma_{i,j}$}
\put(76,31){$\theta_i$}
\put(115,30){$\theta_j$}

\end{picture}
\end{center}

%--------- end of Fig. 1,2------------------------------------------

%\begin{center}
%\begin{picture}(320,140)(-160,-70)
%Fig.1 

%\put(-150,0){\line(1,0){300}}
%\put(40,20){\line(-5,1){100}}
%\put(40,20){\line(-5,-3){100}}
%\put(-60,40){\line(0,-1){80}}

%\put(40,20){\circle*{4}}
%\put(-60,40){\circle*{4}}
%\put(-60,-40){\circle*{4}}
%\put(48,23){$q$}
%\put(-55,47){$p$}
%\put(-55,-50){$p'$}

%\put(97,-40){$H_-$}
%\put(97,40){$H_+$}

%\end{picture}
%\end{center}

\medskip

\begin{cor} {\rm  ({\cite{GFT, Sachs, AKB}})} $\; B(3)=9.$
\end{cor}

\begin{proof} Suppose $P$ is a one-sided kissing arrangement in ${\bf S}^2$ with
$|P|=B(3).$  Since $k(3)=12$, Theorem 1 implies $a+2b\le 12,$ where $a=a(P), b=b(P).$
Note that $A(2,\psi)=\lfloor2\pi/\psi\rfloor.$ Therefore,
$A(2,\arccos(1/\sqrt{3}))=6$.
%$\psi_0:=\arccos(1/\sqrt{3}).$ 
Then Theorem 2 yields
$a\le 6.$ Also we have  $a+b=B(3)\ge 9$. So these nonnegative integer numbers $a$ and $b$ satisfy the following inequalities:
$$a+2b\le 12,\qquad a\le 6, \qquad a+b\ge 9.$$
It is easy to see that there exists only one solution ($a=6,\; b=3$) of these inequalities. Thus,
$B(3)=a+b=9.$
\end{proof}

\begin{cor}
Let $P\subset {\bf S}^3$ be a one-sided kissing arrangement with $|P|=B(4).$ Then there are
only the following possibilities for $|P|, a(P), b(P)$:
$$ |P|=19, \quad (a,b)=(15,4); \; (14,5)$$
$$ |P|=18, \quad (a,b)=(15,3); \; (14,4); \; (13,5); \; (12,6) $$
\end{cor}
\begin{proof} Recently we proved that $k(4)=24$ \cite{Mus,Mus2}. That implies $a+2b\le 24$. On the other hand,
Delsarte's bound gives $A(3,\arccos(1/\sqrt{3}))\le 15$ (see \cite{Boyv}).\footnote{For $\lambda=\arccos(1/\sqrt{3})\approx 54.74^\circ$  Delsarte's bound is stronger than Fejes T\'oth's bound \cite{FeT} $A(3,\lambda)\le 16$. However, Delsarte's bound  also exceeds the lower bound that equals 14. So there is possibility to improve an upper bound for $A(3,\lambda)$ and to skip the case (15,4).}
Therefore, we have the following inequalities for integer numbers $a,b$:
$$ a+2b\le 24,\quad a\le 15, \quad a+b=B(4)\ge 18.$$
There are only the following integer solutions of these inequalities:
$\\(a,b)=(15,4); \; (14,5); \; (15,3); \; (14,4); \; (13,5); \; (12,6).$
\end{proof}

This corollary immediately yields

\begin{cor} {\bf (K. Bezdek)} $\; B(4)=18$ or 19.
\end{cor}

%\noindent{\bf Remark.} Corollary 2 yields that $B(4)=18$ if $(a,b)\ne (15,4), \, (14,5).$ So if we prove it, then we prove our main theorem. 

%The proof of Corollary 2 based on the equality $k(4)=24.$ Actually, for the main theorem
%we don't need it. It is enough to use the Odlyzko - Sloane bound: $k(4)\le 25$ \cite{OdS}. Using this we have just one additional case: $(a,b)=(15,5)$. Since Lemma 2 shows that $(a,b)\ne (15,4)$, we see that $(a,b)\ne (15,5)$.

\section{An extension of Delsarte's method}

For the proof $(a,b)\ne (15,4), (14,5)$ we are using techniques that is called ``Delsarte's method". Here we give a brief introduction to this method (for more details see \cite{Boyv,CS,Del1,Del2,Kab, Lev2,Mus,Mus13,Mus2,OdS,PZ}) and define its extension.

\medskip

Let $P = \{p_1, p_2,\ldots, p_M\}$ be any finite subset of ${\bf S}^{n-1}.$
%\subset{\bf R}^n,\quad {\bf S}^{n-1}=\{x: x\in {\bf R}^n,$ $x\cdot x=||x||^2=1\}.\\$ 
%From here on we will speak of $x\in {\bf S}^{n-1}$ alternatively of points in ${\bf S}^{n-1}$ or of vectors in ${\bf R}^n.$
By $\phi_{i,j}$ we denote the spherical (angular) distance between  $p_i,\,  p_j,$ i.e.
$\phi_{i,j}:=\dist(p_i,p_j).$  

Let us recall the definition of Gegenbauer polynomials by recurrence formula:
$$G_0^{(n)}=1,\quad G_1^{(n)}=t,\quad G_2^{(n)}=\frac{nt^2-1}{n-1},\quad  \ldots,$$ $$G_k^{(n)}=\frac {(2k+n-4)\,t\,G_{k-1}^{(n)}-(k-1)\,G_{k-2}^{(n)}} {k+n-3}$$
Note that $G_k^{(n)}(1)=1.$

\medskip

\medskip

\noindent{\bf Theorem ({ Schoenberg, 1942}).} 
{\it 
If $g_{i,j}=G_k^{(n)}(\cos{\phi_{i,j}}),$ then the matrix $(g_{i,j})$ is positive semidefinite. 

The converse holds also: if $f(t)$ is a real polynomial and for any finite $P\subset{\bf S}^{n-1}$ the matrix $(f(\cos{\phi_{i,j}}))$ is positive semidefinite, then $f$ is a sum of $G_k^{(n)}$ with nonnegative coefficients.}

\medskip

The following theorem is a simple consequence of Schoenberg's theorem.

\medskip

\medskip

\noindent{\bf Theorem {({ Delsarte et al, 1977; Kabatiansky \& Levenshtein, 1978})}.} 
{\em Assume that a polynomial $f$ satisfies the following property:
$$f(t)=f_0G_0^{(n)}(t)+\ldots +f_dG_d^{(n)}(t)\; \mbox{ with }\; f_0\geqslant 0,  \ldots,\, f_d\geqslant 0; \eqno (*)$$
Then for any $P = \{p_1, p_2,\ldots, p_M\}\subset {\bf S}^{n-1}$  we have} 
$$S_f(P):=\sum\limits_{i=1}^M\sum\limits_{j=1}^M{f(\cos{\phi_{i,j}})}\ge f_0M^2.$$ 

\medskip

\begin{theorem} Let $P = \{p_1, p_2,\ldots, p_M\},\; p_i\in{\bf S}^{n-1}$.
Suppose $f$ satisfies $(*)$, and $f(t)\le 0$ for all $t\in [-1,\cos{\psi}]$. Then
$$f_0M^2\le f(1)M+T_f(P),\; \mbox{\; where} \quad T_f(P):=\sum\limits_{(i,j): \phi_{i,j}<\psi,  i\ne j} f(\cos{\phi_{i,j}}),$$
% $$\; I:=\{(i,j): \dist(p_i,p_j)=\phi_{i,j}<\psi, \; i\ne j\}.$$
\end{theorem}
\begin{proof} Since $\phi_{i,i}=0$, we see that $f(\cos{\phi_{i,i}})=f(1)$. 
%Note that $\phi_{i,j}>0$ implies $i\ne j$. 
By the assumption: if 
$\phi_{i,j}\ge \psi$, then
$f(\cos{\phi_{i,j}})\le 0$. Therefore, we have $S_f(P)\le f(1)M+T_f(P).$  On the other hand, $S_f(P)\ge f_0M^2$.
\end{proof}

\begin{cor} {\bf (Delsarte's bound)} Suppose $f$ satisfies the assumptions of Theorem 3 and $f_0>0.$ Then
$$A(n,\varphi)\le \frac{f(1)}{f_0}.$$
\end{cor}
\begin{proof} Let $P$ be a $\varphi$-code of size $M$ on  ${\bf S}^{n-1}$, i.e.
 $\phi_{i,j}\ge\varphi$ for $i\ne j$. Then $T_f(P)=0$. Thus, Theorem 3 yields $f_0M^2\le f(1)M.$ 
\end{proof}

\medskip

\noindent{\bf Remark.} Suitable polynomials $f$ can be found by the linear programming method. See details and algorithms in \cite{CS,Del1,Del2,Kab,Mus2,OdS}.

\medskip

Now we introduce some extensions of Delsarte's method for one-sided kissing arrangements. Note that
all results in this section hold for much more general case: $\varphi$-codes in spherical caps. However, for simplicity, we consider here only one-sided kissing arrangements.% in a hemisphere. 
%Actually, that are a polynomial versions of Theorems 1, 2.

Let $P\subset S_+\subset {\bf S}^{n-1}$ be a one-sided kissing arrangement. For fixed $\theta_0$, 
 where  $60^\circ<\theta_0< 90^\circ$, we define two subsets of $P_a$:
$$ P_{a_1}= \{ p_i\in P_a: \;  \theta_0< \dist(p_i,N)\}, \quad P_{a_2}= \{ p_i\in P_a: \;  \theta_0 \ge \dist(p_i,N)\},$$
$$a_1=a_1(P):=|P_{a_1}|,\quad a_2=a_2(P):=|P_{a_2}|, \quad a=a_1+a_2.$$
 Recall that $|P|=a+b=a_1+a_2+b.$

\begin{theorem} Suppose a polynomial $f$ satisfies $(*)$, $f(t)\le 0$ for all $t\in [-1,1/2]$, and  
$f(t)$ is a monotone increasing function on $[1/2,1]$.
Then
$$(a_1+2a_2+2b)^2f_0\le (a_1+2a_2+2b)f(1)+2a_2f(-\cos{2\theta_0}).$$
\end{theorem}
\begin{proof} Let $$\tilde P=P\bigcup P_{a_2}'\bigcup P_b'.$$ Then
$M=|\tilde P|=|P|+a_2+b=a_1+2a_2+2b.$ 

If $p\in P_{a_2}$, then $\dist(p,p')\ge 180^\circ-2\theta_0.$ For any other pair $(p,q)$ from $\tilde P$ we have  $\dist(p,q)\ge 60^\circ$ (see Theorem 1, Fig. 2). Therefore, the monotonicity assumption implies: $T_f(\tilde P)\le 2a_2f(\cos(180^\circ-2\theta_0))=2a_2f(-\cos{2\theta_0}).$ 
Thus, the proof follows from Theorem 3 with $\psi=60^\circ$.
\end{proof}

\medskip

Now for any function $f(t)$ and $\theta_0, \; 60^\circ <\theta_0<90^\circ$, we introduce  the function $R_f(\theta)$ on the interval $[0,\theta_0]$.

Let $$\psi_0:=\omega(60^\circ,\theta_0,90^\circ)=\arccos\left(\frac{1}{2\sin{\theta_0}}\right). \; $$
For given $e_0\in E\cong {\bf S}^{n-2}$ denote by $C(e_0,\psi_0)$ a spherical cap in the equator $E$ of center $e_0$ and radius $\psi_0$. Let $\Lambda_n(\theta_0)$ be a spherical cone in ${\bf S}^{n-1}$ with vertex $N$ and with the  base $C(e_0,\psi_0)$ (Fig. 4). In other words, $\Lambda_n(\theta_0)$ is a cone of height $90^\circ$ and radius $\psi_0.$

Consider a point $p$ in the height $Ne_0$ with $\dist(N,p)=\theta\le\theta_0.$ Denote by $\Delta_n(\theta,\theta_0)$ the domain in $\Lambda_n(\theta_0)$ that lies below the sphere $S(p)=S(p,60^\circ)$  of center $p$ and radius $60^\circ$. It is easy to see that for $\theta\le\theta_0$ this domain is nonempty.

\medskip

\noindent{\bf Definition.} Consider kissing arrangements $X=\{p_1,\ldots,p_m\}$  in $\Delta_n(\theta,\theta_0)$. Let $\mu(n,\theta,\theta_0)$ be the highest value of $m$. By $\Omega_n(\theta,\theta_0)$ we denote the set of all kissing arrangements in $\Delta_n(\theta,\theta_0).$

% Fig 4
\begin{center}
\begin{picture}(320,140)(-160,-70)

%Fig.4.1 
\put(-95,10){$\theta_0$}
\put(-99,-9){$60^\circ$}
\put(-65,-29){$\psi_0$}
\put(-105,-29){$\psi_0$}
\put(-84,-29){$e_0$}
\put(-84,44){$N$}
\put(-76,-2){$p$}
\put(-77,12){$\theta$}
\put(-95,-48){$\Lambda_n(\theta_0)$}

\put(-80,-20){\circle*{3}}
\put(-80,40){\circle*{3}}
\put(-80,0){\circle*{3}}
\put(-110,-5){\circle*{3}}
%\put(-20,0){\circle*{4}}

%\thicklines
\put(-80,-20){\line(0,1){60}}
%\qbezier (-43,-8) (-80,-32) (-117,-8)
\put(-120,-20){\line(1,0){80}}
\put(-80,-20){\line(-2,1){30}}
%\qbezier (-80,-20) (-97,-18) (-114,-3)
\put(-40,-20){\line(-2,3){40}}
\put(-120,-20){\line(2,3){40}}
%\qbezier (-117,-8) (-110,20) (-80,40)
%\qbezier (-40,-20) (-50,20) (-80,40)

%\put(-60,40){\line(1,-1){40}}
%\put(-140,0){\line(1,1){40}}

%---------------------------
%\put(-95,-65){Fig. 3}
\put(-8,-65){Fig. 4}
%-------------------------
%Fig. 4.2

%\put(40,30){\circle*{4}}
\put(90,-20){\circle*{3}}
\put(90,40){\circle*{3}}
%\put(140,30){\circle*{4}}
%\put(80,20){\circle*{4}}

\qbezier (40,20) (90,-20) (140,20)

%\thicklines
\put(40,-20){\line(1,0){100}}
%\thinlines
\put(40,-20){\line(1,3){20}}
\put(140,-20){\line(-1,3){20}}
%\put(40,30){\vector(1,1){13}}

\put(59,-29){$\psi_0$}
\put(111,-29){$\psi_0$}
\put(86,-29){$e_0$}
\put(94,44){$p$}
\put(18,22){$S(p)$}
\put(72,-12){$\Delta_n(\theta,\theta_0)$}
%\put(86,18){$e_0$}
%\put(25,28){$y_1$}
%\put(147,28){$y_3$}
%\put(50,0){$\Omega_1$}

\end{picture}
\end{center}

Let $Y=\{q_1,\ldots,q_m\}$ be the projection of 
$X=\{p_1,\ldots,p_m\}\subset\Delta_n(\theta,\theta_0)$, where $q_i=\Pi_n(p_i)$, and let 
%denote by $q_i$ the projection of $p_i$ from $N$ onto equator ${\bf S}^{n-2}$. 
$$ F_f(X)=F_f(p_1,\ldots,p_m):= f(\cos{\varphi_1})+\ldots + f(\cos{\varphi_m}), \; 
\varphi_i:=\dist(e_0,q_i).$$
We define $\; R_f(\theta;\theta_0,n)$ as the maximum of  $F_f(X)$ on  $\Omega_n(\theta,\theta_0)$:
$$R_f(\theta)=R_f(\theta;\theta_0,n):=\sup\limits_{X\in \Omega_n(\theta,\theta_0)}{\{F_f(X)\}}.$$
Suppose  $P=\{p_1,\ldots,p_\ell\}\subset{\bf S}^{n-1}$ lies in the  cap of center $N$ and radius $\theta_0.$ Let
$$\tilde R_f(P)=\tilde R_f(P;\theta_0,n):=R_f(\theta_1)+\ldots+R_f(\theta_\ell), \; \, \theta_k:=\dist(N,p_k).$$

\medskip

\noindent{\bf Proposition.} {\em $R_f(\theta)$ is a monotone decreasing function in $\theta.$}
\begin{proof} 
For $\alpha>\beta$, we have $\Delta_n(\alpha,\theta_0)\subset \Delta_n(\beta,\theta_0)$, then $\Omega_n(\alpha,\theta_0)\subset \Omega_n(\beta,\theta_0)$, so then $R_f(\alpha)\le R_f(\beta)$.
\end{proof}

\begin{theorem} Let $\theta_0\in (60^\circ,90^\circ)$. Suppose a polynomial $f$ satisfies the following assumptions:
$f=f_0G_0^{(n-1)}+\ldots +f_dG_d^{(n-1)}\; \mbox{ with }\; f_0\ge 0,  \ldots,\, f_d\ge 0; $ and \\
$f(t)\le 0 \; $  for all $\; t\in [-1,\cos{\psi_0}].$\\ 
%$f(t) \;$  is a monotone increasing function on $\; [\cos{\psi_0},1].$\\
Then for any one-sided kissing arrangement $P\subset {\bf S}^{n-1}$ with $|P|=M$  we have
$$f_0M^2\le Mf(1)+2a_2R_f(60^\circ;\theta_0,n)+2\tilde R_f(P_b;\theta_0,n).$$
\end{theorem}
\begin{proof}
First let us prove the theorem for the case when  the North pole $N$ is not a point in $P$. 
%any other $\tilde p\in P$ we have $\dist(N,\tilde p)=\dist(p,\tilde p)\ge 60^\circ$. This implies that $\tilde p\in P_a$ and $|P_b|=1$. It contradicts the assumption:  $|P_b|=5$. 
In this case the projection $Q$ of $P$ is well defined. 
%Let $Q$ be the projection of $P.$
%Denote this projection by $Q=Q_a\bigcup Q_b$. 

If we apply Theorem 3 for $f$ and $Q=\Pi_n(P)\subset E$, then we obtain
$$f_0M^2\le S_f(Q)\le f(1)M+T_f(Q).$$

From the assumption: $f(t)\le 0, \; t\in [-1,\cos{\psi_0}]$, for  $q_i, q_j\in Q$ it follows that 
if $\dist(q_i,q_j)\ge \psi_0$, then $f(\cos{\phi_{i,j}})\le 0$. 
Arguing as in Theorem 2  it is easy to see that for any distinct $q_1, q_2\in Q_{a_1}$ we have
$\dist(q_1,q_2)\ge \psi_0$.
%$:=\omega(60^\circ,71.74^\circ,90^\circ)\approx 58.2299^\circ.$$

For $q_1,q_2\in Q_{a_2}$ we have
$$
\dist(q_1,q_2)\ge \omega(60^\circ,\theta_0,\theta_0)>\omega(60^\circ,\theta_0,90^\circ)=\psi_0.
$$
In the case $q_1, q_2\in Q_b$ we have 
$$
\dist(q_1,q_2)>\omega(60^\circ,60^\circ,60^\circ)>\psi_0.
$$ 
%(see Theorem 3 in \cite {Mus2}). =\arccos(1/3)\approx  70.5288^\circ

Therefore, $\dist(q_1,q_2) < \psi_0$ only if $$
(i) \; \, q_1\in Q_{a_1}, \; q_2\in Q_{a_2},\; \mbox{ or } \; \, (ii) \; q_1\in Q_a, \; q_2\in Q_b.
$$

Let $p\in P_{a_2}, \; \dist(N,p)=\theta, \; e_0=\Pi_n(p)$, and let $X=\{p_1,\ldots,p_m\}\subset P$ with $\dist(e_0,q_i)< \psi_0$. 
By $(i)$ we have $X\subset P_{a_1}$. This implies that $X$ lies below $S(p,60^\circ)$. Thus, $X\subset\Delta_n(\theta,\theta_0)$ and $F_f(X)\le R_f(\theta)$.
The same holds for $p\in P_b$ (the case $(ii)$).
% \; X\subset P_a$. In this case, no more points from $P_b$ can lie in $\Lambda_n(\theta,\theta_0).$ 
We obtain
$$T_f(Q)\le 2\tilde R_f(P_{a_2})+2\tilde R_f(P_b).$$
 By definition $\dist(N,p)\ge 60^\circ$ for $p\in P_{a_2}$. Since $R_f(\theta)$ is a monotone decreasing function,
we see that $\tilde R_f(P_{a_2})\le a_2R_f(60^\circ)$.
%Consider a point $q\in Q_{a_2}$ and two points $p,r\in Q_{a_1}$ such that $\dist(q,p)<\psi_0,$ and $\dist(q,r)<\psi_0 $. It is not hard to prove that $\angle pqr>\omega(60^\circ,\psi_0,\psi_0)\approx 72.0474^\circ$. Since $\angle pqr>72^\circ$, there are at most four points from $Q_{a_1}$ in the spherical cap of center $q$ and radius $\psi_0$.

It can be shown in the usual way that the theorem for $N\notin P$ implies the theorem for   $N=p\in P$. Indeed, if we shift $p$  to some point $p(\varepsilon)$ which is a small distance $\varepsilon$ away  from $N$, then we can apply the theorem for a  $\varphi$-code  
$P(\varepsilon)=p(\varepsilon)\cup P\setminus p$, where $\varphi=60^\circ-\varepsilon.$ Since $R_f(\theta)$ is a continuous
function and $R_f(0)$ is well defined the theorem for $N\in P$ follows  whenever $\varepsilon\to 0.$
\end{proof}

In the end of this section we consider some geometric properties of kissing arrangements in $\Delta_n(\theta,\theta_0)$.
 These properties will help us to compute $R_f(\theta)$
in the next section.

Note that for $\theta> 30^\circ$  the sphere $S(p,60^\circ)$ intersects the equator. That yields
 $$
\dist(e_0,q_i)\ge \rho(\theta):= \omega(60^\circ,\theta,90^\circ)=\arccos{\left(\frac{1}{2\sin{\theta}}\right)}, \quad \theta> 30^\circ
$$

It is easy to see that $p_i\in\Delta_n(\theta,\theta_0)$ implies
$\dist(N,p_i)\ge \gamma(\theta,\theta_0)$, where the function $\gamma=\gamma(\theta,\theta_0)=\gamma_1(\theta,60^\circ,\psi_0)$ can be defined by the equation
$$\omega(60^\circ,\theta,\gamma)=\psi_0, \quad \gamma\ge\theta.$$
In other words, $\gamma_1(\theta,\varphi,\psi)$ is the length of  $AB$ of a spherical triangle $ABC$ with $|AC|=\theta,\; |BC|=\varphi,$ and $\angle BAC=\psi$. In fact, this equation has two solutions. Let us denote the second solution by $\gamma_2(\theta,\varphi,\psi)$.

We get
$$\dist(q_i,q_j)\ge  d(\theta)=d(\theta,\theta_0):=\omega(60^\circ,\gamma(\theta,\theta_0),90^\circ), \quad i\ne j.$$

Therefore, the set $Y=\{q_i\}$ is a $d(\theta)$-code on the cap $C(e_0,\psi_0)\subset E\cong{\bf S}^{n-2}.$ Thus,
$$\mu(n,\theta,\theta_0)\le A(n-1, d(\theta,\theta_0),\psi_0).$$

It was proved (see Theorem 3 in \cite {Mus2}) that for $\varphi>\psi$ we have
$$A(n,\varphi,\psi)=A(n-1,\omega(\varphi,\psi,\psi)).$$
By $\theta_*(\theta_0)$ denote the solution of the equation $ d(\theta,\theta_0)=\psi_0$. We obtain
$$\mu(n,\theta,\theta_0)\le A(n-2,\omega( d(\theta),\psi_0,\psi_0))   \; \mbox{ for } \; \theta>\theta_*(\theta_0).$$

Consider the case $n=4$. We have a $d(\theta)$-code  $Y=\{q_1,\ldots,q_m\}$  on the cap $C(e_0,\psi_0)\subset {\bf S}^{2}.$ Now we consider two cases: a) $\theta>\theta_*(\theta_0)$ and 
b) $\theta\le\theta_*(\theta_0)$.

\medskip

\noindent a) By definition for $\theta>\theta_*(\theta_0)$ we have $d(\theta)>\psi_0$. It can be shown  (see Corollary 2 in \cite{Mus2}) that in this case $m\le 5$.

Let $m=5$. Then $q_i$ are vertices of a spherical convex pentagon \cite[Theorem 2]{Mus2}.  We denote this pentagon by $V_5$.  Since $\dist(q_i,q_j)\ge d(\theta)$ we have that the length of any side of $V_5$ is at least $d(\theta)$. On the other hand, $V_5\subset C(e_0,\psi_0).$ 
%Let us to find the lower bound on $\dist(e_0,q_i)$
It is not hard to prove that the minimum  distance between $e_0$ and a vertex $q_i$ of $V_5$ achieves  when $V_5$ is an equilateral pentagon of sides lengths $d(\theta)$  and four other vertices $q_j$ lie in the boundary of the cap $C(e_0,\psi_0)$. 
That yields
$$
\dist(e_0,q_i)\ge \rho_5(\theta,\theta_0):=\gamma_2(\psi_0,d(\theta), \eta),\quad
 \eta=180^\circ-3\omega(d(\theta),\psi_0,\psi_0)/2.
$$

\noindent b) In this case $m\le 6$. The maximal size arrangement can be achieved when $q_6=e_0$ and $\{q_1,\ldots,q_5\}$ are vertices of a spherical convex pentagon.

Indeed, let  
$\delta=\omega(d(\theta),\psi_0,\psi_0).$ It can be shown numerically that for any $\theta_0\in(60^\circ,90^\circ)$ we have $\delta>60^\circ$. Let $V$ be the (spherical) convex hull of $Y$. It can be shown that $\angle q_ie_0q_j\ge \delta$, where $q_i, q_j$ are vertices of $V$.
  Since $\lfloor 360^\circ/\delta\rfloor\le 5$ we have that the number of vertices of $V$ is at most 5. It's not hard to prove that no more than one $q_i\in Y$ can lie inside $V$. Thus, $m\le 6.$

\medskip

These facts can be summarized as follows:

\begin{theorem} Suppose $X=\{p_1,\ldots,p_m\}$ is a kissing arrangment in  $\Delta_n(\theta,\theta_0)$. Let $q_i=\Pi_n(p_i)$. Then

\medskip

\noindent $
1) \; \dist(q_i,q_j)\ge  d(\theta,\theta_0), \; \mbox{ and }\; \dist(e_0,q_i)\ge \rho(\theta) \; 
\mbox{ for }\; \theta>30^\circ;
$  

\medskip

\noindent$2) \; \, m\le \mu(n,\theta,\theta_0)\le A(n-1,d(\theta,\theta_0),\psi_0);$ 

\medskip

\noindent$3) \; \, m\le \mu(n,\theta,\theta_0)\le A(n-2,\omega(d(\theta,\theta_0),\psi_0,\psi_0))   \; \mbox{ for } \; \theta>\theta_*(\theta_0);$

\medskip

\noindent$4)$ if $n=4$, then $m\le\mu(4,\theta,\theta_0)\leq 6$; moreover
for $\theta>\theta_*(\theta_0)$ we have $m\le 5,$  and if $m=5$, then
$\; \dist(e_0,q_i)\ge \rho_5(\theta,\theta_0).$

\end{theorem}

\section{B(4)=18}

In this section we prove the main theorem of this paper:
\begin{theorem} $\; B(4)=18.$
\end{theorem}

Now we prove that $(a,b) \ne (15,4), (14,5).\; $ In this section $\theta_0=71.74^\circ.$
So we have  $\psi_0\approx 58.2299^\circ$, $\theta_*=\theta_*(\theta_0)\approx 25.8526^\circ,$ and $\cos{2\theta_0}=-\cos{36.52^\circ}$.\\
\noindent{\em Here and below numbers are shown to four decimal places.}

%$$ P_{a_1}= \{ p_i\in P_a: \;  71.74^\circ < \dist(p_i,N)\}, \quad P_{a_2}= \{ p_i\in P_a: \;  71.74^\circ \ge \dist(p_i,N)\}.$$
%Let $a_1(P):=|P_{a_1}|,\; a_2(P):=|P_{a_2}|$. Then $a_1(P)+a_2(P)=a(P).$
%$$a_2(P)=|P_2|,\quad a_1(P)=a(P)-a_2(P),\quad (a_1,a_2,4)$$

\begin{lemma} Let $P\subset {\bf S}^3$ be a one-sided kissing arrangement with $|P|=19.$\\
$1)$ If $(a,b)=(15,4)$, then $a_2\le 5.$\\
$2)$ If $(a,b)=(14,5)$, then $a_2\le 3.$
\end{lemma}
\begin{proof} 
Let $f_{OS}(t)=f_0+f_1G_1^{(4)}(t)+\ldots+f_9G_9^{(4)}(t),\; $  where
$\\ f_0=1, \; f_1=3.6181,\; f_2=6.1156,\; f_3=7.0393,\; f_4=5.0199,\\ 
f_5=2.313, \; f_6=f_7=f_8=0, \; f_9=0.4525.$

%\medskip

%\medskip

The polynomial $f_{OS}$ satisfies the assumptions of Theorem 4.\footnote{This polynomial was applied by Odlyzko and Sloane \cite{OdS} to prove that $k(4)\le 25.$ Indeed, Delsarte's bound gives $k(4)=A(4,60^\circ)\le  f_{OS}(1)/f_0=f_{OS}(1)\approx 25.5584.$ Then $\, k(4)\le 25.$} Now we apply this theorem for the proof. 
%\medskip
%Note that $K=2a_2.$ Suppose $M=|\tilde P|=29,\;$ then $a_2\le 6,\; K\le 12.$ From Corollary 5 
%On the other hand, in the sum $S_f(Q)$ we have only $|Q|+2a_2$ positive terms. Namely, $|Q|$ terms 
%$f(\cos{\phi_{i,i}})=f(1)$, and $2a_2$ terms $f(\langle p,p'\rangle)$, where $p\in P_{a_2}$. For all other we have $\cos{\phi_{i,j}}\le 1/2$, so  1) yields $f(\cos{\phi_{i,j}})\le 0$. From $a_2\le 6$

\noindent 1) Assume the converse. Then $\; a_2\ge 6,$ and
$M=a_1+2a_2+2b=|P|+a_2+b\ge 29.$ 
For $a_2=6$ from Theorem 4 
$$841=f_029^2 \le 29f_{OS}(1)+12f_{OS}(\cos(36.52^\circ))\approx 840.8819$$
follows, a contradiction.

\noindent 2) If $(a,b)=(14,5)$ and $a_2=4$, then $M=28$. We have
$$784=f_028^2 \le 28f_{OS}(1)+8f_{OS}(\cos(36.52^\circ))\approx 782.0941,$$
a contradiction.

Since $f_0M^2-Mf_{OS}(1)-2a_2f_{OS}(\cos(36.52^\circ))$ is increasing whenever $M$ is increasing, we complete our proof.
\end{proof}

\begin {lemma} $\; (a,b) \ne (15,4).$
\end {lemma}
\begin{proof} Assume the converse. 
%Then $|P_a|=15.$
%Let $Q_a=Q_{a_1}\bigcup Q_{a_2}$ be the projection of $P_a=P_{a_1}\bigcup P_{a_2}$ onto equator ${\bf S}^2.$ 
%From this follows that if $K$ is the number of all pairs $(p,q), p, q\in Q_a$ with $\dist(p,q)<\psi_0$, then $K\le8a_2.$ (Here $(p,q)$ and $(q,p)$ are two different pairs.)
%From Lemma 1 follows that $a_2\le 5$. Thus, $K\le 40$.

Let $g(t)=g_0+g_1G_1^{(3)}(t)+\ldots+g_9G_9^{(3)}(t),\; $  where
$\\ g_0=1, \; g_1=2.7986 ,\; g_2=3.6388,\; g_3=3.4429,\; g_4=2.1227, \;
g_5=0.8637, \\ g_6=g_7=g_8=0, \; g_9=0.1281$. Then $g(1)=g_0+g_1+\ldots+g_9= 13.9948.$

The polynomial $g$ satisfies the assumptions of Theorem 5 with $n=4$. Note that $g(t)$ is a monotone increasing function on $[\cos{\psi_0},1]$.

%,\; \psi=\psi_0,\; \theta=\arccos(1/\sqrt{3}).\; $ Note that $g(1/\sqrt{3})\approx 0.3668.$ Then 
Let us apply Theorem 5 for $P_a$. Since $|P_a|=15$ and $b(P_a)=0$ we have
$$225=g_015^2\le 15g(1)+2a_2R_g(60^\circ).$$

Now we prove that $R_g(60^\circ)\le 4g(1/\sqrt{3}).$ 

Theorem 6 yields $$\dist(e_0,q_i)\ge \rho(60^\circ)=\arccos(1/\sqrt{3})\approx 54.7356^\circ,$$
$\mu(4,60^\circ,\theta_0)\le 5$, and for $m=5$ we have
$$\dist(e_0,q_i)\ge\rho_5=\rho_5(60^\circ,\theta_0)\approx  57.5043 ^\circ.$$  
Then the monotonicity of $g(t)$ on $[\cos{\psi_0},1]$ implies for $m=5$: $F_g(X)\le 5g(\cos{\rho_5})$ and for $m\le 4$: $F_g(X)\le 4g(1/\sqrt{3}).$
%For $m=5$ it is not hard to prove that $\dist(q_0,q_i)>\alpha_0= 57.5043 ^\circ.$
Since $$4g(1/\sqrt{3})\approx 1.4673>5g(\cos{\rho_5})\approx 0.332,$$ we have $R_g(60^\circ)\le 4g(1/\sqrt{3}).$  (It is easy to see that $R_g(60^\circ)= 4g(1/\sqrt{3}).$)

By Lemma 1, $a_2\le 5$, so that
$$225\le 15g(1)+2a_2R_g(60^\circ)\le 15g(1)+40g(1/\sqrt{3})\approx 224.5946$$
 - a contradiction.
\end{proof}

\begin{lemma} $\; (a,b) \ne (14,5).$\footnote{There is old conjecture that the maximal kissing arrangement in four dimensions is  unique up to isometry, i.e. is the ``24-cell"  $P_{24}$.  From this conjecture easily follows Lemma 3. Indeed, let $\hat P:=P\bigcup P'_b$, then $|\hat P|=a+2b=24.$ The conjecture yields that $\conv(\hat P)$ is the ``24-cell". Since at most 18 vertices of $P_{24}$ can lie in a closed hemisphere, $|P|\le 18$ - a contradiction.} %  However, the conjecture is not proven yet.}
\end{lemma}
\begin{proof}
Assume the converse. 

\noindent {\bf 1.} Let $P_b=\{p_1,\ldots,p_5\},\; \theta_i:=\dist(p_i,N)$ , and
%Without of generality it can be assumed that 
let  $\{\bar p_i\}$ be sorted $\{p_i\}$ by $\theta_i:$ 
$$\{p_i\}=\{\bar p_j\}, \quad \bar\theta_i:=\dist(\bar p_i,N), \quad 60^\circ>\bar\theta_1\ge\bar\theta_2\ge\bar\theta_3\ge\bar\theta_4\ge\bar\theta_5>0.$$ 
By $P_b(k)$ denote $\{\bar p_1,\ldots, \bar p_k\}$. Here $1\le k\le 5$, and  $P_b(5)=P_b.$
%Suppose  is the projection of $Z_k=\{z_1,\ldots, z_k\}\subset P_b$  onto ${\bf S}^2$ from $N$.  
%Consider any $k$ points from $Q_b$. Let us denote this subset of $Q_b$ by $Y_k$. 
Let $$P(k):=P_a\bigcup P_b(k), \quad Q(k):=Q_a\bigcup Q_b(k).$$ 
Then $|P(k)|= |Q(k)|=14+k.$

Suppose a polynomial $f$ satisfies the assumptions of Theorem 5 with $n=4$. 
%and $f(t)$ is a monotone increasing function on $[\cos{\psi_0},1]$.
%satisfies $(*)$ with $n=3$, and\\ (i) $f(t)\le 0$ for all $t\in [-1,\cos{\psi_0}]$, \\
%(ii) $f(t)$ is a monotone increasing function on $[\cos{\psi_0},1].$\\ 
%(Here as above $\psi_0=\omega(60^\circ,71.74^\circ,90^\circ)\approx 58.2299^\circ.$)

Then Theorem 5 for $P(k)$ yields  
$$f_0(14+k)^2\le S_f(Q(k))\le (14+k)f(1)+2a_2R_f(60^\circ)+2\tilde R_f(P_b(k)).\eqno (1)$$

Let $h(t)=h_0+h_1G_1^{(3)}(t)+\ldots+h_9G_9^{(3)}(t),\; $  where\\
$ h_0=1, \; h_1=2.722 ,\; h_2=3.6191,\; h_3=3.4323,\; h_4=2.4227,\;
h_5=1.1991,$\\  $h_6=h_7=h_8=0, \; \; h_9=0.0195. \quad \; h(1)=h_0+h_1+\ldots+h_9= 14.7413.$

Note that  $h$ (as well as $g$) satisfies the assumptions of Theorem 5, and $h(t)$ is a monotone increasing function on $[\cos{\psi_0},1]$. Actually, $g$ and $h$ are close to each other.
As above, $R_h(60^\circ)=4h(1/\sqrt{3})$. However, for $k\le 3$ using $f=g$ in $(1)$ we get better results than for $f=h$, but $h$ is better than $g$ for $k=4,5$. 
 
%Therefore, we can apply $(2)$ for $g$ and $h$. For $f=g$ and $f=h$ we can give an explicite expressions for $R_f(\theta).$

%By definition $T_f(\tilde Q(k))=\sum f(\cos{\phi_{i,j}})$, where $0<\phi_{i,j}<\psi_0.$
%Let $K$ be the number of all pairs $(\tilde q_i,\tilde q_j)$ from $Q_a$ with $\phi_{i,j}<\psi_0$. Lemma 1 yields $a_2\le 4.$ Then $K\le 8a_2\le 32.$
%For any distinct $\tilde q_i,\tilde q_j\in Q_a$ we have $\phi_{i,j}\ge\arccos(1/\sqrt{3}).$ So (ii) implies $f(\cos{\phi_{i,j}})\le f(1/\sqrt{3}).$ 

%Let
%$$\tilde R_f(Q_b(k),Q_a):=\sum {f(\cos{\phi_{i,j}})}, \mbox{ where } \; q_i\in Q_b(k),\; \tilde q_j%\in Q_a, \; 0<\phi_{i,j}<\psi_0.$$
%Then $T_f(\tilde Q(k))\le 2\tilde R_f(Q_b(k),Q_a)+32f(1/\sqrt{3})$.
%We obtain
By Lemma 1 $a_2\le 3$. Therefore, for $f=g$ and $f=h$ we have
$$ (14+k)^2\le (14+k)f(1)+24f(1/\sqrt{3})+2\tilde R_f(P_b(k)). \eqno (2) $$ 

Recall that
$\tilde R_f(P_b(k))= R_f(\bar\theta_1)+\ldots+R_f(\bar\theta_k).$
Thus, $(2)$ yields
$$R_f(\bar\theta_1)+\ldots+R_f(\bar\theta_k)\ge k^2/2+14k+98-(7+k/2)f(1)-12f(1/\sqrt{3}). \eqno (3)$$

\medskip

\noindent {\bf 2.} Now we compute the upper bound on $R_f(\theta)$ for $f=g,h$, and
$0\le\theta\le 60^\circ.$

Let  $\theta\in[0^\circ,60^\circ]$.  If $X$ is a kissing arrangement in $\Delta_4(\theta,\theta_0)$, then Theorem 6 yields some constraints for $Y$, where $Y$ is the projection of $X$. Using these constraints we introduce the function $L_f(\theta)$ such that $R_f(\theta)\le L_f(\theta).$

 For a fixed $e_0\in {\bf S}^2$ we consider the set $\Psi(\theta)$ of all $d(\theta)$-codes $Y=\{q_1,\ldots,q_m\}$ on ${\bf S}^2$ such that
$m\le 6; \;
 \varphi_i:=\dist(e_0,q_i)\ge \rho(\theta)$ for $\theta>30^\circ;\; $  $m\le 5$ for $\theta>\theta_*$ and if $m=5$ then $\varphi_i\ge \rho_5(\theta).$ 
%Denote by $\Psi$ the set of all $Y\subset {\bf S}^2$ that satisfy these constraints. 
Let
$$ H_f(Y):= \tilde f(\varphi_1)+\ldots + \tilde f(\varphi_m),\quad 
L_f(\theta):=\sup\limits_{Y\in\Psi(\theta)}{\{H_f(Y)\}},$$
where
$$ \tilde f(\varphi):=\left\{
\begin{array}{l}
f(\cos{\varphi}) \quad  0\le\varphi\le\psi_0\\
0 \qquad \quad \; \;  \varphi>\psi_0
\end{array} 
\right.
$$
Theorem 6 implies that if $X\in \Omega_4(\theta,\theta_0)$, then $Y\in\Psi(\theta)$. Thus $R_f(\theta)\le L_f(\theta)$.

Consider some $\theta\in [30^\circ,60^\circ]$.  Recall that $ f(\cos{\varphi})$ is a monotone decreasing function on $[0,\psi_0]$.    Then  (as well as in Lemma 2)  we have for
$m=5$: $H_f(Y)\le 5\tilde f(\rho_5(\theta))$, and for $m\le 4: H_f(Y)\le 4\tilde f(\rho(\theta)).$ It can be shown numerically that in this interval $4\tilde f(\rho(\theta))>5\tilde f(\rho_5(\theta))$. Therefore,
 $L_f(\theta)\le 4\tilde f(\rho(\theta)).$ 

 Let  
$\delta_4\approx  44.4577^\circ$  is defined by the equality: $\cos^2{\rho(\delta_4)}=\cos{d(\delta_4)}$. This equality means that a  spherical square of side length $d(\delta_4)$ has the circumradius  equals $\rho(\delta_4)$. Then  for $\theta\in [\delta_4,60^\circ]$ the side length of a  spherical square with the circumcircle of center $e_0$ and radius $\rho(\theta)$ is not less than $d(\theta)$. Thus,
$$
L_f(\theta)=4\tilde f(\rho(\theta))=4f\left(\frac{1}{2\sin{\theta}}\right)\; \mbox{ for } \; \theta\in [\delta_4,60^\circ].
$$

If $\theta<\delta_4$, then the circumradius of  a  square of side length $d(\theta)$  is greater than $\rho(\theta)$. 
%We say that a kissing arrangement $X\subset\Delta(\theta,\theta_0)$ is optimal if   
It have been proved (see Theorem 4 \cite{Mus2}) that if $Y$ is optimal (i.e. $H_f(Y)=L_f(\theta)$)  and $|Y|=4$, then   $Y$   is a spherical rhomb $q_1q_2q_3q_4$ with edge lengths $d(\theta)$. 
%this inequality implies the fact that 
%$|X|\le 4$.
% Therefore,
%$
%R_f(\theta)=\tilde f(\varphi_1)+\ldots + \tilde f(\varphi_4),
%$

% is the distance between $e_0$ and the vertex $q_k$  of $Y$.   

Consider a rhomb $Y=q_1q_2q_3q_4$ of side lengths $d(\theta)$ with $\dist(q_1,q_3)\ge\dist(q_2,q_4)$ such that the circumradius of $q_1q_2q_4$ is $\rho(\theta)$. Under these conditions this rhomb there exists and uniquely defined for $\theta\in [\delta_3,\delta_4].$ Here $\delta_3\approx 37.4367^\circ$ is defined by the condition that 
$\rho(\delta_3)$ is the circumradius of a regular triangle of side lengths $d(\delta_3)$. It can be verified numerically  that for all $\theta\in [\delta_3,\delta_4]\; $      $H_f(Y)$ achieves its maximum when $e_0$ is the circumcenter of $q_1q_2q_4$.  Since for this  arrangement  $\dist(e_0,q_3)$  is uniquely defined by $\theta$ we denote it by $r_4(\theta)$. Then 
$$
L_f(\theta)=3f\left(\frac{1}{2\sin{\theta}}\right) + \tilde f(r_4(\theta)) \; \mbox{ for } \; \theta\in [\delta_3,\delta_4].
$$

Note that for $\theta=\delta_3$ we have $\varphi_3>\psi_0$, i.e.  $\tilde f(\varphi_3)=0.$ Computations show that for $\theta\in [30^\circ,\delta_3]$ the maximum of  $H_f(Y)$ is achieved when $q_1q_2q_4$ is a regular triangle  of side lengths $d(\theta)$ with $\dist(e_0,q_1)=\rho(\theta)$. Let
$r_3(\theta):=\dist(e_0,q_2)=\dist(e_0,q_4).$ Since in this case $\varphi_3>\psi_0$, we obtain
%\footnote {Actually, we have one more constraint for $Y$: if $|Y|=3$ ( i.e. $Y$ is a triangle), then there are $i,j$ such that $\dist(q_i,q_j)>d(\theta)$. Therefore for $\theta\in (30^\circ,\delta_3)$ we have $R_f(\theta)<L_f(\theta)$.}
$$
 L_f(\theta)=f\left(\frac{1}{2\sin{\theta}}\right) + 2\tilde f(r_3(\theta)) \; \mbox{ for } \; \theta\in [30^\circ,\delta_3].
$$

For $\theta=30^\circ$ we see that $q_1=e_0$ and $r_3(\theta)=d(\theta)> \psi_0$. Therefore,
$L_f(30^\circ)=f(1).$  Note that for $\theta\in [30^\circ, \theta_*]$ we have $f(1)>5\tilde f(\rho_5(\theta))$. Then, as above, we can assume that $Y$ is a rhomb. 
Computations show that  $H_f(Y)$ attains its maximum when $e_0$ is a vertex of $Y$, i.e. $L_f(\theta)=f(1)$. For $\theta<\theta_*$  it is not hard to show that $H_f(Y)$ achieves its maximum when  $e_0\in Y$ and other $q_i$ lie at the distance $d(\theta)$ away from $e_0$.
Thus,
$$
L_f(\theta)=f(1)+5\tilde f(d(\theta)) \; \mbox{ for } \; \theta\in [0,30^\circ].
$$
(Note that if $\theta\ge\theta_*$, then $d(\theta)\ge\psi_0$, i.e. $\tilde f(d(\theta))=0.$)

%$H(Q_a,Q_b):=\sum\limits_{(p,q)\in\Lambda} h(\langle p,q\rangle), \quad 
%\Lambda=\{(p,q): p\in Q_a,\, q\in Q_b,\; \dist(p,q)<\psi_0\}. $$ Here by $\langle p,q\rangle$ denoted the inner product of unit vectors $p$ and $q$, i.e. $\langle p,q\rangle=\cos(\dist(p,q))$.   In other words, $H(Q_a,Q_b)$ is the sum of all positive terms   in $S_h(\tilde Q)$ with $q_i\in Q_a, \; q_j\in Q_b.$

%\medskip

%\noindent {\bf Claim.} $H(Q_a,Q_b)\le 34.9039$

%\medskip

%Note that from this Claim follows Lemma 3. Indeed,
% We have $h(1/\sqrt{3})\approx  0.2632$. Thus,
%$$361\le 19h(1)+2H(Q_a,Q_b)+32h(1/\sqrt{3})<  358.32,$$
%a contradiction.

\medskip

\noindent {\bf 3.} Consider the domain $\Theta=:\{(\theta_1,\ldots,\theta_5)\}$ in ${\bf R}^5.$ The proof of Lemma 3 follows from the fact that  $\Theta$ is empty set under the constraints $(3)$ and $\dist(p_i,p_j)\ge 60^\circ.$

It is not so hard to find bounds for $\theta_k$ under the constraints $(3)$. For instance,
since $R_f(\theta)$ is a monotone decreasing function, we have
$$
R_f(\bar\theta_1)+\ldots+R_f(\bar\theta_k)\le kR_f(\bar\theta_k)\le kL_f(\bar\theta_k).
$$
Combining this with $(3)$ we get
$$
\bar\theta_1 < 53.93^\circ, \;  \bar\theta_2 < 48.45^\circ, \; \bar\theta_3 < 46.72^\circ, \; 
\bar\theta_4 < 45.53^\circ, \; \bar\theta_5 < 44.47^\circ.
$$
Here for $k=1,2,3$ we applied $f=g$, and for $k=4,5$: $f=h.$

Moreover, using convexity of $L_f(\theta)$ on $[\delta_4,60^\circ]$ we can give inequalities for sums of $\theta_k$. Perhaps, it is possible  using these ineaqualities to prove the lemma without computer works. However, we think that is a hard problem.  Here for the proof we apply the numerical method
that was developed in our paper \cite{Mus2}.

  Denote by $P_5(\alpha)$ a convex polytope in ${\bf S}^3$ with five vertices (see Fig. 5). The lengths of all edges $P_5(\alpha)$,  except $p_2p_4,\; p_3p_5$, are equal to $60^\circ,$ and $$\dist(p_2,p_4)=\alpha\in [60^\circ,90^\circ],
\; \dist(p_3,p_5)=\beta,  \; 2\cos{\alpha}\cos{\beta}+\cos{\alpha}+\cos{\beta}=0.
$$
%(\theta_1,\ldots,\theta_5)=
%and $R(P_b):=\rho_h(\theta_1)+\ldots+\rho_h(\theta_5).$ Then $H(Q_a,Q_b)\le R(P_b).$
%Actually, for $R$ we have completely the same optimization problem that is considered in \cite{Mus2}.  For this problem an optimal $P_b$ belongs to the 1-parametric family in  (Fig. 4).  Then our method (see Section 5 in \cite{Mus2}) gives $\max\{R(P_b)\}\approx 34.9039.$
%By assumption $\theta_i\le 60^\circ.$ (Recall that $\dist(y_i,y_j)\ge 60^\circ.$
\begin{center}
\begin{picture}(320,140)(-80,-70)
% Fig. 
\put(57,-65){Fig. 5: $P_5(\alpha)$}

\put(10,-20){\circle*{5}}

\put(90,-40){\circle*{5}}
\put(70,20){\circle*{5}}

\put(130,40){\circle*{5}}

\put(70,60){\circle*{5}}

\thicklines
\put(10,-20){\line(4,-1){80}}
\put(10,-20){\line(3,2){60}}

\put(70,20){\line(3,1){60}}
\put(90,-40){\line(1,2){40}}

\put(70,60){\line(-3,-4){60}}
\put(70,60){\line(0,-1){40}}
\put(70,60){\line(1,-5){20}}
\put(70,60){\line(3,-1){60}}

\thinlines
\multiput(90,-40)(-1,3){20}%
{\circle*{1}}

\put(69,-14){$\alpha$}

%labels
\put(-5,-21){$p_5$}
\put(97,-41){$p_2$}
\put(57,24){$p_4$}
\put(134,44){$p_3$}
\put(73,65){$p_1$}
 
\end{picture}
\end{center}

Let $X=\{p_1, \ldots, p_5\}$  be a kissing arrangement in the cap $C(N,60^\circ)\subset{\bf S}^3$, and 
let  $$F(X)=F_1(\theta_1)+\ldots+F_5(\theta_5), \; \theta_i=\dist(N,p_i),$$
where $F_i(\theta)$ is a monotone decreasing function in $\theta$. 
From Theorem 5 \cite{Mus2} it follows that $F$  attains its maximum on $\{X\}$ at $X_0$ when  $\conv(X_0)$ (the convex hull of $X_0$) is isometric to $P_5(\alpha)$ for some  $\alpha\in [60^\circ,90^\circ]$. Note that Theorem 5 is not be assumed a strong monotonicity of $F_i(\theta)$.

Let us apply this fact with $F(P_b)=L_f(\bar\theta_1)+\ldots+L_f(\bar\theta_k)$, i.e.
$F_i=L_f$ for $i=1,\ldots,k$, and $F_i=0$ for $i>k.$ (It's easy to see that $L_f(\theta)$ is a monotone decreasing function.)
Thus, we can assume that $\conv(P_b)=P_5(\alpha).$

Let us show how using this fact to obtain a contradiction for the lemma. 

Let the vertices $p_1, p_2, p_3$ of $P_5(\alpha$) be fixed. Then the vertices $p_4=p_4(\alpha)$ and
 $p_5=p_5(\alpha)$ are uniquely determined by $\alpha$. Consider the following domain
$$
D_\varepsilon(u_1,u_2,u_3):= \{p\in S^3: \; 
\theta_i=\dist(p,p_i)
\in [u_i,u_i+\varepsilon]\}, \quad
\varepsilon>0.
%\theta_2\in [v,v+\varepsilon], \; \theta_3\in [w,w+\varepsilon]\}%
$$
In general position 
$D=D_\varepsilon(u_1,u_2,u_3)$ looks like a curvilinear parallelepiped, i.e $D$ has (at most) eight vertices. 

Let us denote by $S=S(u_1,u_2,u_3;\alpha)=\{s_i\}$ the set that consists of vertices of $D$ and points on the boundary of $D$ which are intersections of the great circles through $p_ip_j, i\ne j, j=4,5$, 
and $D$. For small $\varepsilon$ and $u_i>0$ it is easy to see that $|S|\le 10$.
Lemma 5 \cite{Mus2} yields that the function $\theta_i=\dist(p,p_i),\; i=4,5,$ in $p$ achieves its minimum on $D$ at some $s_j$. 
Let us denote this point by $\tilde p_i(\alpha,D)$.
%(Actually, for $i=1,2,3$ it follows from the definition: $\theta_1\ge u,\; \theta_2\ge v, \; \theta_3\ge w$.) 
Thus, $L_f(\theta_i)$ (as well as any monotone decreasing function in $\theta$) achieves its maximum on $D$ at $\theta_i=\tilde\theta_i(\alpha):=\dist(p_i,\tilde p_i(\alpha,D))$. 

Note that $\theta_4$ is increasing and $\theta_5$ is decreasing whenever $\alpha$ is increasing. Then
for $\alpha\in[r,r+\varepsilon]$ we have $\tilde\theta_4(\alpha)\ge \tilde\theta_4(r)$ and 
$\tilde\theta_5(\alpha)\ge \tilde\theta_5(r+\varepsilon)$. Finally, for $L_f(\theta_i)$ on
$\bar D_\varepsilon=D\times[r,r+\varepsilon]$ we have
$$
L_f(\theta_i)\le L_f(u_i),\; i=1,2,3; \; L_f(\theta_4)\le L_f(\tilde\theta_4(r)); \; 
L_f(\theta_5)\le L_f(\tilde\theta_5(r+\varepsilon))\eqno (4).
$$

Using $(4)$ we can easily to find  an upper bound on $L_f(\theta_i)$ on any $\bar D_\varepsilon$. Then we can check $(3)$ for points in $\bar D_\varepsilon$. So if we cover $S_+\times[60^\circ,90^\circ]$ by $\bar D_\varepsilon$ we can check $(3)$ for all possible $p, \; \{p_i\}$ and $\alpha$.
Our computations show that there no points in $\Theta$ that satisfy $(3)$ - a contradiction.
\end{proof}

\section {Conclusions}

It is clear that between kissing numbers and one-sided kissing numbers there are some relations.
Look at these nice equalities:

$$ n=2, \quad 4=B(2)=\frac{k(1)+k(2)}{2}=\frac{2+6}{2};$$
$$ n=3, \quad 9=B(3)=\frac{k(2)+k(3)}{2}=\frac{6+12}{2};$$
$$ n=4, \quad 18=B(4)=\frac{k(3)+k(4)}{2}=\frac{12+24}{2}.$$

We are not expect that the same holds for all dimensions, i.e.
$$B(n)=\bar K(n):=\frac{k(n-1)+k(n)}{2}$$
for all $n.\; $ However, there are many reasons that $B(n)=\bar K(n)$ for
$ n=5, 8, 24;\; \;$ $B(n)\approx \bar K(n)$ for all $n$; and this equality holds asymptotically
$$ \lim_{n\to\infty}{\frac{B(n)}{\bar K(n)}}=1.$$

If we consider the minimal vectors of $D_5, E_8,$ and Leech lattices, then we see that $B(5)\ge 32, \; B(8)\ge 183, \; B(24)\ge 144855$. We propose the following 

\medskip

\medskip

\noindent{\bf Conjecture.} $B(5)=32, \; B(8)=183, \; B(24)=144855.$

\medskip

\medskip

Note that  our method could be applied to higher dimensions. However, it does not give sharp upper bound on $B(n).$ It's an interesting problem to find better method, in particular for $n=5, 8, 24.$

\medskip

\medskip

\medskip

\medskip

\noindent{\bf Acknowledgment.} I wish to thank K\'aroly Bezdek and G\'abor Fejes T\'oth for helpful discussions and useful comments on this paper.

\end{document}